\documentclass[12pt, twoside, a4paper]{article}

\ifx\TestIngCommand\undefined\relax
\else\input ../../../proc/ppreamb-book.tex 
\input init.tex \fi
\let\TestIngCommand\undefined

\usepackage{amsmath,amscd,amssymb} 
\usepackage{graphics}
\newtheorem{theo}{Theorem}
\newtheorem{lem}{Lemma}

\newskip\ttglue\ttglue=.5em plus.25em minus.15em
\def\firstname#1{\def\FIRSTNAME{#1}\ignorespaces}
\def\lastname#1{\def\LASTNAME{#1}\ignorespaces}           
\def\middleinitial#1{\def\MIDDLEINI{#1}\ignorespaces}
\def\department#1{\def\DEPARTMENT{#1}\ignorespaces}
\def\institute#1{\def\INSTITUTE{#1}\ignorespaces}
\def\address#1{\def\ADDRESS{#1}\ignorespaces}      
\def\country#1{\def\COUNTRY{#1}\ignorespaces}                   
\def\otheraffiliation#1{\def\OTHERAFFILIATION{#1}\ignorespaces}
\def\email#1{\def\EMAIL{#1}\ignorespaces}
\newcount\autcount\autcount=0
\newcount\affcount\affcount=0
\newcount\numcount\newcount\nummcount
\newcount\nummmcount\newcount\nummmmcount                                  
\def\writename#1#2{\ \kern-1ex\hbox{                                     
  \csname AUthor\the#1\endcsname\                                          
  \edef\TESTSTR{}\expandafter\ifx\csname auTHor\the#1\endcsname\TESTSTR
  \else\csname auTHor\the#1\endcsname.\ \fi                                
  \csname authOR\the#1\endcsname$^{\csname AFF\the#1\endcsname}$   
  \expandafter\ifx\csname corr\number#1\endcsname\relax 
  \else\thanks{Corresponding author.}\ \fi                                
  }\ifnum#1<#2, \else\ \kern-1ex\fi}
\def\writeemail#1{
  \nummcount=0\relax\nummmcount=0\relax
  \loop\ifnum\nummcount<\autcount\advance\nummcount by1\relax
    {\expandafter\ifnum\csname AFF\the\nummcount\endcsname=#1\relax
    \global\advance\nummmcount by1\fi}\repeat
  \nummcount=0\relax\nummmmcount=0\relax
  \loop\ifnum\nummcount<\autcount\advance\nummcount by1\relax
    {\expandafter\ifnum\csname AFF\the\nummcount\endcsname=#1\relax
    \global\advance\nummmmcount by1\relax\def\blank{}\expandafter
    \ifx\csname EMAIL\the\nummcount\endcsname\blank(no e-mail)
    \else\csname EMAIL\the\nummcount\endcsname
    \fi                                
    \ifnum\nummmmcount<\nummmcount; \fi\fi}\repeat}
\long\def\BeginAuthorList#1\EndAuthorList{#1\relax
  \author{\vbox{\hsize=390pt\noindent\numcount=0\relax
    \loop\ifnum\numcount<\autcount\advance\numcount by1\relax
      \writename{\numcount}{\autcount}
      \repeat}\\[2mm]
    \vbox{\small\numcount=0\relax
      \loop\ifnum\numcount<\affcount\advance\numcount by1\relax
        \vbox{{\count0=\numcount\relax
          \loop\expandafter\ifnum\csname AFF\the\count0\endcsname
            <\numcount\relax\advance\count0 by1\relax\repeat
          $^{\csname AFF\the\count0\endcsname}$}
        \def\BLANK{}\expandafter\ifx\csname DEPT\the\numcount\endcsname 
          \BLANK                      
          \else\csname DEPT\the\numcount\endcsname, \fi 
        \csname INST\the\numcount\endcsname, 
        \csname ADDR\the\numcount\endcsname,
        \csname COUN\the\numcount\endcsname
        \edef\TEST{}\expandafter\ifx\csname OTHE\the\numcount\endcsname
          \TEST                                                            
          .\else;\break\csname OTHE\the\numcount\endcsname.\fi}
        \vbox{\writeemail{\numcount}}
        \repeat}\\}}
\expandafter\def\csname x1\endcsname{}
\expandafter\def\csname x2\endcsname{}
\expandafter\def\csname x3\endcsname{}
\expandafter\def\csname x4\endcsname{}
\expandafter\def\csname x5\endcsname{}
\expandafter\def\csname x6\endcsname{}
\expandafter\def\csname x7\endcsname{}
\expandafter\def\csname x8\endcsname{} 
\expandafter\def\csname x9\endcsname{}
\def\Author#1#2{\global\advance\autcount by1\relax#2
  \expandafter\edef\csname AUthor\the\autcount\endcsname{\FIRSTNAME}
  \expandafter\edef\csname auTHor\the\autcount\endcsname{\MIDDLEINI}
  \expandafter\edef\csname authOR\the\autcount\endcsname{\LASTNAME}
  \expandafter\edef\csname EMAIL\the\autcount\endcsname{\EMAIL}
  \let\tempera\"\def\"{\string\"}\expandafter\ifx\csname x\DEPARTMENT
    \endcsname\relax
    \global\advance\affcount by1\relax\let\"\tempera
    \expandafter\edef\csname DEPT\the\affcount\endcsname{\DEPARTMENT}
    \expandafter\edef\csname INST\the\affcount\endcsname{\INSTITUTE}
    \expandafter\edef\csname ADDR\the\affcount\endcsname{\ADDRESS}
    \expandafter\edef\csname COUN\the\affcount\endcsname{\COUNTRY}
    \expandafter\edef\csname OTHE\the\affcount\endcsname{\OTHERAFFILIATION}
    \expandafter\edef\csname AFF\the\autcount\endcsname{\the\affcount}
  \else\expandafter\edef\csname AFF\the\autcount\endcsname{\DEPARTMENT}
  \fi\let\"\tempera\ignorespaces}
\def\CorrespondingAuthor#1#2{
  \expandafter\xdef\csname corr\number#1\endcsname{cor}
  \Author#1{#2}}
\def\PaperTitle#1{\title{\bf#1}}
\def\Category#1{\ignorespaces}
\def\keywords#1{{\noindent \emph{Keywords:}
  \def\BLANK{}\def\TEST{#1}\ifx\BLANK\TEST(n/a).\else#1\fi}}
\begin{document}
\PaperTitle{Note on Spectral Factorization Results of Krein and Levin}
\Category{(Pure) Mathematics}
\date{}
\BeginAuthorList
  \Author1{
    \firstname{Wayne}
    \lastname{Lawton}
    \middleinitial{M}
    \department{Department of the Theory of Functions, Institute of Mathematics and Computer Science}
    \institute{Siberian Federal University}
    \otheraffiliation{}
    \address{Krasnoyarsk}
    \country{Russian Federation}
    \email{wlawton50@gmail.com}}
\EndAuthorList
\maketitle
\begin{abstract}
Bohr proved that a uniformly almost periodic function $f$ has a bounded spectrum if and only if it extends to an entire function $F$ of exponential type $\tau(F) < \infty$. If $f \geq 0$ then a result of Krein implies that $f$ admits a factorization $f = |s|^2$ where $s$ extends to an entire function $S$ of exponential type $\tau(S) = \tau(F)/2$ having no zeros in the open upper half plane. The spectral factor $s$ is unique up to a multiplicative factor having modulus $1.$ Krein and Levin constructed $f$ such that $s$ is not uniformly almost periodic and proved that if $f \geq m > 0$ has absolutely converging Fourier series then $s$ is uniformly almost periodic and has absolutely converging Fourier series. We derive neccesary and sufficient conditions on $f \geq m > 0$ for $s$ to be uniformly almost periodic, we construct an $f \geq m > 0$ with non absolutely converging Fourier series such that $s$ is uniformly almost periodic, and we suggest research questions.
\end{abstract}
\noindent{\bf 2010 Mathematics Subject Classification:47A68;42A75;30D15}
\footnote{\thanks{This work is supported by the Krasnoyarsk Mathematical 
Center and financed by the Ministry of Science and Higher Education 
of the Russian Federation in the framework of the establishment and 
development of Regional Centers for Mathematics Research and 
Education (Agreement No. 075-02-2020-1534/1).}}
\section{Notation}
$:=$ means `is defined to equal' and iff means `if and only if'.
$\mathbb N = \{1,2,3,...\},$ $\mathbb Z,$ $\mathbb Q,$ $\mathbb R,$ $\mathbb C$
are the natural, integer, rational, real, and complex numbers. For $z \in \mathbb C,$
$x := \Re z; y := \Im z$ are its real; imaginary coordinates.
$\mathbb D := \{z \in \mathbb C: |z| \leq 1\}$ is the closed unit disk,
$\mathbb D^o := \{z \in \mathbb C: |z| \leq 1\}$ is the open unit disk,
and
$\mathbb T := \{z \in \mathbb C: |z| = 1\}$ is the circle. 
$\mathbb U := \{ \, z \in \mathbb C\, : \, \Im z \geq 0 \, \}$ is the closed upper half--plane and $\mathbb U^o := \{ \, z \in \mathbb C\, : \, \Im z \geq 0 \, \}$ is the open upper half plane. For $z \in \mathbb C,$ $\chi_z(x) := e^{izx}.$ For closed $K \subset \mathbb C,$ $C_b(K)$ is the $C^*$--algebra of bounded continuous complex--valued functions on $K$ with norm 
$||f|| := \sup_{z \in K} |f(z)|.$ For $\rho > 0,$
$\mathfrak{D}_\rho : C_b(\mathbb R) \rightarrow C_b(\mathbb R)$ is the dilation operator $(\mathfrak{D}_\rho f)(x) := f(\rho \, x).$
Zeros of nonzero entire functions are denoted by sequences $z_n, n \geq 1$ of finite (possibly zero) or infinite length. 
\\ \\
Sections 2 and 3 records classical results and derives preliminary results. Section 4 presents two new results. Section 5 suggests areas for research.
\section{Entire Functions of Exponential Type}
\begin{lem}\label{herglotz}
If $E : \mathbb U \rightarrow \mathbb C$ is continuous, holomorphic on $\mathbb U^o,$
$\Re \, E$ is bounded above, and $\Re \, E|_{\mathbb R} = 0,$
then $E(z) = i(a+bz)$ for $a$ real and $b \geq 0.$
\end{lem}
{\bf Proof} For some $c > 0$ the function $E_1 : \mathbb D \backslash \{-1\} \rightarrow \mathbb R$ defined by
$$
E_1(z) := c - \Re E\left( \frac{i-iz}{1+z} \right), \ \ z \in \mathbb D\backslash \{-1\}
$$
is positive and continuous and its restriction to $\mathbb D^o$ is harmonic. Herglotz \cite{herglotz}, (\cite{rudin},Theorem 11.30) proved that there exists a unique positive Borel measure $\mu$ on $\mathbb T$ such that
$$
    E_1(z) := \Re \int_{w \in \mathbb T} \frac{w + z}{w-z} \, d\mu(w), \ \ z \in \mathbb D^o.
$$
Since $\Re \, E|_{\mathbb R} = 0$ implies $E_1(z) = c$ for $z \in \mathbb T \backslash \{-1\},$
it follows that there exists $b \geq  0$ with
$d\mu(w) = cd\sigma + b \delta_{-1}$ where $\sigma$ is normalized Haar measure on $\mathbb T$ and $\delta_{-1}$ is the point measure at $-1.$ Therefore
$$
    E_1(z) = c + b\,  \Re \, \frac{1-z}{1+z}, \ \ z \in \mathbb \mathbb D \backslash \{-1\} \implies
  \Re E(z) = \Re i b z, \ \ z \in \mathbb U
$$
and the result follows from the Cauchy--Riemann equations.
\\ \\
Throughout this section $F$ denotes a nonconstant entire function. Its order (\cite{boas}, Definition 2.1.2), defined by
$$
\rho(F) := \limsup_{r \rightarrow \infty} \frac{\log \log M(r)}{\log r}
$$
where 
$$
M(r) := \max \{ \, |F(z)| \, : \, |z| \leq r \, \}, \ \ 0 \leq r,
$$
is in $[0,\infty].$
For integers $k \geq 1$ the entire function $z^k;e^{z^k}$ has order $0;k.$
The convergence exponent (\cite{boas}, Definition 2.5.2) $\rho_1 \in [0,\infty]$ of the zeroes $z_j, j \geq 1$ of $F$ is
$$
	\rho_1 := 
	\inf \left \{ 
	\alpha \in (0,\infty] : \sum_{n} |z_n|^{\alpha} < \infty
	\right \}.
$$ 
If $\rho_1 < \infty$ the genus $p$ of its zeros is the smallest nonnegative integer with
$$
	\sum_{n} |z_n|^{p+1} < \infty.
$$ 
Clearly $p \in [\rho_1,\rho_1-1].$
Finite sequences have $\rho_1 = p = 0.$
The infinite sequence $\frac{1}{n};\frac{1}{n \log^2(n+1)}$ has $\rho_1 = 1;1$ and $p = 1;0.$
\\ \\
The following result shows how to construct certain entire functions with finite order (\cite{boas}, (2.6.3), Theorem 2.6.4).
\begin{lem}\label{canprod}
If $z_j, j \geq 1$ has finite convergence exponent $\rho_1$ and genus $p,$ then the canonical product of genus $p$ defined by
$$P(z) := \prod_{0 < |z_n|} E(z/z_n,p)$$
where 
$$E(z,p) := (1-z)\exp \left[ z + \frac{z^2}{2} + \cdots + \frac{z^p}{p} \right]$$
is the Weierstrass primary factor, converges uniformly on compact subsets to an entire function $P$ that has order $\rho_1$ and zeros $\{z_j: 0 < |z_n|\}.$
\end{lem}
Hadamard's factorization (\cite{boas}, Theorem 2.7.1) gives
\begin{lem}\label{Hadamard}
If $F$ has finite order $\rho(F)$ then its zeros $z_j, j \geq 0$ have finite convergence exponent $\rho_1 \leq \rho(F).$ If $m \geq 0$ is the multiplicity of $0$ as a root of $F,$ then there exists a polynomial $Q$ having degree $q \leq \rho(F)$ such that
$$
F(z) = z^m e^{Q(z)} P(z)
$$ 
where $P$ is the canonical product of genus $p.$ Furthermore $\rho(F) = \max \{q, \rho_1\}.$
\end{lem}
If $F$ has finite order $\rho$ we define its type 
$$\tau(F) := \limsup_{r \rightarrow \infty} r^{-\rho} \log M(r)$$
and say it has finite type if $\tau(F) < \infty.$
Define
$$n_F(r) := \hbox{cardinality of } \{j : |z_j| \leq r\}.$$
Lindel\"{o}f \cite{lindelof}, (\cite{boas}, Theorem 2.10.1) proved
\begin{lem}\label{lindelof}
If $F$ is an entire function whose order $\rho$ is a positive integer, then $F$ has finite type iff both $n_F(r) = O(r^\rho)$ and
$$
	\sup_{r \geq 0} \left| \sum_{0 < |z_n| \leq r} z_{n}^{-\rho} \right| < \infty.
$$
\end{lem}
We say $F$ has exponential type if $\rho(F) = 1$ and $\tau(F) < \infty.$
For $\alpha \in \mathbb C \backslash \{0\},$ the function $e^{\alpha z}$ has exponential type $\tau = |\alpha|.$
\\ \\
Krein \cite{krein} proved
\begin{lem}\label{krein}
	If $F$ is of exponential type and $f := F|_{\mathbb R}$ is bounded and nonegative, then $f$ admits a factorization $f = |s|^2$ where $s$ extends to an entire function $S$ of exponential type that has no zeros in $\mathbb U^{o}.$ Moreover $\tau(S) = \tau(F)/2$ and $s$ is unique up to multiplication by a constant having modulus $1.$
\end{lem}
Levin (\cite{levin}, p. 437) said that Krein used approximation of $f$ by Levitan trigonometric polynomials (\cite{levin}, Appendix I, Section 4). We observe that $F = S\overline S$ where $\overline S(z) := \overline {S(\overline z)}.$
Levin (\cite{levin}, Chapter V) defines
$$
	A := \left \{ F \hbox{ entire } : 
	\sum_{0 < |z_n|} \left| \, \Im \frac{1}{z_n}\, \right| < \infty \right \}.
$$
Boas (\cite{boas}, p. 134) proved:
\begin{lem}\label{boas}
If $F$ is of exponential type, then $F \in A$ iff
$$
	\sup_{L \geq 1} \int_{1}^L x^{-2}|\log |f(x)f(-x)|\, dx f(x)\, \}\, dx < \infty.
$$
Moreover, this condition holds whenever
$$
	\int_{-\infty}^{\infty} (1+x^2)^{-1} \max \{0, \log |f(x)|\} dx < \infty.
$$
\end{lem}
Ahiezer \cite{ahiezer} extended Krein's result by proving:
\begin{lem}\label{ahiezer}
If $f = F|_{\mathbb R} \geq 0,$ then $F \in A$ iff $f$ admits a factorization $f = |s|^2$ where
$s$ extends to an entire function $S$ of exponential type $\tau(S) = \tau(F)/2$ having no zeros in $\mathbb U^{o}.$
\end{lem}
We summarize Levin's proof in (\cite{levin}, p. 437). The if part follows from Lemma \ref{lindelof}. The roots of $F$ occur in conjugate pairs so we order $\{z_n\}$ so $n$ odd $\implies \Im n_n \leq 0$ and $z_{n+1} = \overline z_n.$
Lemma \ref{Hadamard} implies there exists an integer $m \geq 0$ and $a, b \in \mathbb R$ such that
$$
	F(z) = z^{2m}\, e^{2az+2b}\,  \prod_{n \geq 1} \left(1 - \frac{z}{z_n}\right)\, e^{\frac{z}{z_n}}.
$$
Define $\gamma := -\sum_{n \ odd} \Im \frac{1}{z_n},$ and $s := S|_{\mathbb R}$ where
$$
	S(z) :=  z^m\, e^{az+b+i\gamma z} \, \prod_{n \ odd} \left(1 - \frac{z}{z_n}\right)\, e^{\frac{z}{z_n}}.
$$
Then $F = S \overline S$ and $S$ has no zeros in $\mathbb U^{o}.$ Since $S \in A,$ Lemma \ref{lindelof} implies that $S$ is of exponential type. The rest of the proof shows that
$\tau(S) = \tau(F)/2.$
\\ \\
$\xi(\mathbb R)$ is the space of smooth complex valued functions on $\mathbb R$ with the topology of uniform convergence of derivatives of every order on compact subsets. Its dual space
$\xi^{\, \prime}(\mathbb R)$  is the space of compactly supported distributions (\cite{hormander}, Theorem 1.5.2). For $u \in \xi^{\, \prime}(\mathbb R)$ we define $\alpha(u), \beta(u) \in \mathbb R$ so that $[\alpha(u),\beta(u)]$ is the smallest closed interval containing the support of $u,$ and we define its Fourier-Laplace transform by $\widehat u(z) := u(\chi_{-z}), z \in \mathbb C.$
\begin{lem}\label{PWS}
Conditions 1 and 2 are equivalent and imply condition 3
\begin{enumerate}
\item $F$ is of exponential type and $\exists \, c > 0, N \in \mathbb Z$ with
$
|F(x)| \leq c\, (1 + |x|)^N, \ x \in \mathbb R.
$
\item $F = \widehat u$ for $u \in \xi^{\, \prime}(\mathbb R)$ with
$\tau(F) = \max \{ |\alpha(u)|, |\beta(u)| \}.$
\item $\lim_{r \rightarrow \infty} r^{-1} n_F(r) =\frac{\beta(u)-\alpha(u)}{2\pi}.$
\end{enumerate}
\end{lem}
Schwartz  \cite{schwartz} proved the equivalence of conditions 1 and 2 and Titchmarsh
\cite{titchmarsh} proved they imply condition 3.
\section{Uniformly Almost Periodic Functions}
$T(\mathbb R)$ is the algebra of trigonometric polynomials spanned by
$\{\chi_\omega :\omega \in \mathbb R\},$
and the
$C^*$--algebra $U(\mathbb R)$
of uniformly almost periodic functions is its closure with respect to $|| \cdot ||.$
Bohr \cite{bohr1} proved that if $f \in U(\mathbb R)$ its mean value
$$
	M(f) := \lim_{L \rightarrow \infty} \, \frac{1}{2L} \int_{-L}^L f(x)dx
$$
exists, defined its Fourier transform $\widehat f : \mathbb R \rightarrow \mathbb C$ by
$$
	\widehat f(\omega) := M(f\,\chi_{-\omega}), \ \omega \in \mathbb R,
$$
and proved that its spectrum
$
	\Omega(f) := \{w \in \mathbb R : \widehat f \neq 0\}
$
is countable. For $f \in U(\mathbb R)$ its Fourier series is the formal sum
$$
f \sim \sum_{\omega \in \Omega(f)} \widehat f(\omega) \, \chi_{\omega},
$$
its bandwidth
$
	b(f) := \sup \Omega(f) - \inf \Omega(f) \in [0,\infty],
$
and
$$
	||f||_A := \sum_{\omega \in \Omega(f)} |\widehat f(\omega)| \in [0,\infty].
$$
We define the following subsets of $U(\mathbb R):$
\begin{enumerate}
\item Bandlimited algebra $B(\mathbb R) := \{\,f \in U(\mathbb R) \, : \,
b(f) < \infty \, \}.$
\item Wiener Banach algebra $A(\mathbb R) := \{f \in U(\mathbb R) \, : \,
||f||_A < \infty\}.$
\item Hardy Banach algebra $H(\mathbb R) := \{f \in U(\mathbb R) \, : \,
\Omega(f) \subset [0,\infty) \}.$
\item Invertible Hardy functions $IH(\mathbb R) := \{f \in H(\mathbb R) : 
\frac{1}{f} \in H(\mathbb R) \}.$
\end{enumerate}
Bohr \cite{bohr1}, (\cite{levin}, Chapter VI, p. 268, Corollary to Theorem 1) proved:
\begin{lem}\label{bohr1}
	If $f \in U(\mathbb R)$ then $f \in B(\mathbb R)$ iff $f$ extends to an entire function $F$ of exponential type. Then $\tau(F) = \max \{|\inf \Omega(f)|, |\sup \Omega(f)|\}.$
\end{lem}
The Fourier series for $f \in U(\mathbb R)$ converges absolutely iff $f \in A(\mathbb R).$
Cameron \cite{cameron} and Pitt \cite{pitt}, (\cite{grs}, Section 29, Corollary 2 to Theorem 2) proved
\begin{lem}\label{cameronpitt}
	If $f \in A(\mathbb R)$ and $\Phi$ is holomorphic in an open region containing the closure of
	$f(\mathbb R),$ then the composition $\Phi \circ f \in A(\mathbb R).$
\end{lem}
If $f \in B(\mathbb R)$ and $f \geq 0$ then Lemma \ref{bohr1} implies that $f$ extends to an entire function $F$ of exponential type $\tau(F) = b(f)/2$ so Lemma \ref{krein} implies that $f = |s|^2$ where $s$ extends to an entire function $S$ of exponential type $\tau(S) = b(f)/4.$
Levin (\cite{levin}, Appendix 2, Theorem 2) used Lemma \ref{cameronpitt} to prove
\begin{lem}\label{levin1}
If $f \in A(\mathbb R) \cap B(\mathbb R)$ and there exists $m > 0$ with $f \geq m,$ then the spectral factor $s \in A(\mathbb R) \cap B(\mathbb R).$
\end{lem}
For $\Delta > 0$ let $[\Delta]$ denote the set of entire functions $F$ of exponential type $\tau(F) = \Delta$ such that $f := F|_{\mathbb R} \in B(\mathbb R)$ and $-\Delta, \Delta \in \Omega(f).$
Krein and Levin obtained a precise characterization of the zeros of functions in $[\Delta]$ and published these results without proofs in \cite{kreinlevin}. In (\cite{levin}, Appendix VI) for the first time they gave proofs for these results and used them (\cite{levin}, Appendix VI, p. 463) to prove
\begin{lem}\label{kreinlevin}
There exists $\Delta > 0$ and $F \in [\Delta]$ with $f \geq 0$ whose spectral factor $s \notin U(\mathbb R).$
\end{lem}
Levin's result (\cite{levin}, Chapter V1, Section 2, Lemma 3) implies:
\begin{lem}\label{asym}
If $h \in U(\mathbb R)$ then $-\infty < \Delta := \inf \Omega(h)$ iff $h$ extends to a continuous function $H$ on $\mathbb U$ which is holomorphic on $\mathbb U^o$ and satisfies
$$
	\lim_{y \rightarrow \infty} e^{-i \Delta (x+iy) } \, H(x+iy) = \widehat h(-\Delta)
$$
where convergence is uniform in  $x.$ Therefore $h \in H(\mathbb R)$ iff $H \in C_b(\mathbb U).$ 
\end{lem}
The Poisson kernel functions $P_y \, : \, \mathbb R \rightarrow \mathbb R, y > 0$ are
$$
	P_y(x) := \frac{y}{\pi} \frac{1}{x^2 + y^2}, \  \ \ x + iy \in \mathbb U^o.
$$
For $f \in C_b(\mathbb R)$ its Poisson integral $P[f] : \mathbb U \rightarrow \mathbb C$ is
$P[f](x) := f(x)$ and
$$
    P[f](x+iy) := \begin{cases}
				f(x), \  \  \ y = 0 \\
			\int_{-\infty}^{\infty} P_y(x-s)\, f(s) \,  ds, \  \  \ y > 0
			\end{cases}
$$
\begin{lem}\label{poisson}
If $f \in C_b(\mathbb R),$ then $P[f] \in C_b(\mathbb U),$ 
its restriction $P[f]|_{\mathbb U^o}$ is harmonic, and
$$
	\sup_{x \in \mathbb R} \Re f(x) \geq \Re P[f](z) \geq \inf_{x \in \mathbb R} \Re f(x), \ \ z \in \mathbb U^o.
$$
If $f \in U(\mathbb R)$ then $f \in H(\mathbb R)$ iff $P[f]|_{\mathbb U^o}$ is holomorphic.
\end{lem}
{\bf Proof}
The first assertion follows since $P_y(x)$ is harmonic, positive valued, and $\int P_y(x)dx = 1, y > 0.$
The second assertion follows since
$$
	P[\chi_\omega](z) = \begin{cases} e^{i\omega z}, \  \  \  \omega \geq 0 \\
	e^{i\omega \overline z}, \  \  \  \omega < 0. \end{cases}
$$
Bohr \cite{bohr1}, (\cite{levin}, Chapter VI, Theorem 2) proved:
\begin{lem}\label{bohr2}
	If $h \in H(\mathbb R)$ is nonzero, $H := P[h],$ and $z_n, n \geq 1$ are the zeros of $H,$ then $\{\Im z_n\}$ is bounded iff $\inf \Omega(h) \in \Omega(h).$
\end{lem}
Bohr \cite{bohr2}, (\cite{levin}, p. 274, footnote) proved:
\begin{lem}\label{bohr3} If $h \in U(\mathbb R)$ and $|h|^2 \geq m$ for some $m > 0,$
then there exists $c \in \mathbb R$ and $\theta \in U(\mathbb R)$ such that
$$
(\arg h)(x) = cx + \theta(x), \  \ x \in \mathbb R.
$$
\end{lem}
\begin{lem}\label{invh1}
If $h \in IH(\mathbb R)$ then $|P[h]|$ is bounded below by a positive number, 
$\Re \log P[h] \in C_b(\mathbb U),$ and $\widehat h(0) \neq 0.$
\end{lem}
{\bf Proof} Lemma \ref{asym} implies
$P[h], P[1/h], P[h]P[1/h] \in C_b(\mathbb U).$ 
Lemma \ref{poisson} implies
$P[h]P[1/h]$ is holomorphic on $\mathbb U^o.$ Since
$P[h]P[1/h]|_{\mathbb R} = 1,$ 
the Schwarz reflection principle (\cite{rudin}, Theorem 11.14) 
implies $P[h]P[1/h] = 1.$ 
Therefore $|P[h]|$ is bounded below by a positive number, so $\log P[h]$
exists, is unique up to addition by an integer multiple of $2\pi i,$ and
$\Re \log P[h] = \log |P[h]| \in C_b(\mathbb U).$
Since $P[h]$ and $P[1/h]$ have no zeros,
Lemma \ref{bohr2} implies $\inf \Omega(h) \in \Omega(h)$
and $\inf \Omega(1/h) \in \Omega(1/h).$ Since $\{0\} = \Omega(1) \subset \Omega(h) + \Omega(1/h)$ it follows that $\inf \Omega(h) = \inf \Omega(1/h) = 0$ so $\widehat h(0) \neq 0.$
\begin{lem}\label{invh2}
If $h \in IH(\mathbb R),$ $f := |h|^2,$ and $f \in B(\mathbb R),$ 
then $h \in B(\mathbb R)$ and $\chi_{-b(f)/4}\, h$ is a spectral factor of $f.$
\end{lem}
{\bf Proof}
Lemma \ref{krein} implies $f = |s|^2$ where the spectral factor $s$ extends to an entire function $S$ of exponential type
$\tau(S) = b(f)/4$ which has no zeros in $\mathbb U^0.$
Define $S_1(z) := e^{ib(f)z/4}\, S(z)|_{\mathbb U}$
and $H := P[h].$
Then $S_1 \in C_b(\mathbb U)$ is bounded and holomorphic 
with no zeros in ${\mathbb U^o}$ and Lemma \ref{invh1} implies that
$H$ is holomorphic on $\mathbb U^o$ and $|H|$ is bounded below by a positive number. Therefore $G := S_1/H \in C_b(\mathbb U)$ is  holomorphic with no zeros on
$\mathbb U^o,$ and $|G(x)| = |s(x)|/|h(x)| = 1$ for $x$ real. Therefore $E := \log G$ exists and satisfies the hypothesis of Lemma \ref{herglotz} hence $E(z) = i(a+bz)$ for some $a$ real and $b \geq 0.$ Therefore
$
	e^{ib(f)z/4}\, S(z) = e^{i(a+bz)}\, H(z)
$
hence
$
	s = e^{ia} \chi_{b - b(f)/4}\, h
$
so $h \in B(\mathbb R).$ Lemma \ref{invh1} implies  $\inf \Omega(h) = 0.$
Since $\inf \Omega(s) = -b(f)/4$ it follows that $b = 0$ hence $\chi_{-b(f)/4}\, h$ is a spectral factor. 
\\ \\
Boas (\cite{boas}, Theorem 11.1.2) proved this generalization of Sergei Bernstein's classic theorem \cite{bernstein}, (\cite{boas}, Theorem 11.1.1) for polynomials:
\medskip
\begin{lem}\label{bernstein} If $F$ is an entire function of exponential type 
and $f := F|_{\mathbb R}$ is bounded, then 
$f^{\, \prime} := \frac{df}{dx}$ is bounded and 
$||f^{\, \prime}|| \leq \tau(F) \, ||f||.$
\end{lem}
\begin{lem}\label{cauchyriemann}
If $f \in B(\mathbb R),$ $m > 0,$ $f \geq m,$ 
$g := \frac{1}{2} \log f,$ $u := P[g],$ $u_o := u|_{\mathbb U^o},$ and
$v_o : \mathbb U^o \rightarrow \mathbb R$ is a harmonic function conjugate to $u_o,$ 
then $v_o$ is uniformly continuous on each horizontal line and on each vertical ray in 
$\mathbb U^o$ so extends to a continuous function 
$v : \mathbb U \rightarrow \mathbb R.$ Furthermore $|v(z)| = O(|z|).$
\end{lem}
{\bf Proof} Let $\gamma := \frac{1}{\pi}\int_{-\infty}^{\infty} \frac{|s^2-1|}{(s^2+1)^2}ds.$ $\frac{\partial v_o}{\partial x}$ is bounded on each horizontan line since 
$$
	\left|\frac{\partial v_o}{\partial x}(x+iy)\right| = \left|\frac{\partial u_o}{\partial y}(x+iy)\right| \leq \frac{\gamma}{2y}
	\max \{|\log m|,\, |\log ||f||\, | \} < \infty, \  \  \ x+iy \in \mathbb U^o.
$$
$\frac{\partial v_o}{\partial y}$ is bounded on each vertical ray since Lemma \ref{bernstein} implies
$$
\left|\frac{\partial v_o}{\partial y}(x+iy)\right| = \left|\frac{\partial u_o}{\partial x}(x+iy)\right|
\leq \frac{b(f)\, ||f||}{4m} < \infty, \  \  \ x+iy \in \mathbb U^o.
$$
\section{New Results}
\begin{theo}\label{result1}
If $f,$ $g,$ $u$ and $v$ are as in Lemma \ref{cauchyriemann},  
then $f$ has a spectral factor $s \in B(\mathbb R)$ iff 
$v|_{\mathbb R}$ has the form in Lemma \ref{bohr3}.
\end{theo}
{\bf Proof}
Let $H := e^{u+iv},$ $h := H|_{\mathbb R}.$ If $v|_{\mathbb R} = \arg h$ has the form in Lemma \ref{bohr3} then
$
h = \chi_{c}\, \sqrt f \, e^{i\theta} \in U(\mathbb R).
$ 
Since $P[h] = H$ is holomorphic on $\mathbb U^o,$ 
Lemma \ref{poisson} implies $h \in H(\mathbb R).$ 
A similar argument implies $\frac{1}{h} \in H(\mathbb R)$ so $h \in IH(\mathbb R).$ 
Since $|h|^2 = f,$ Lemma \ref{invh2} implies $h \in B(\mathbb R)$ and 
$s := \chi_{-b(f)/4}\, h$ is a spectral factor of $f$ in $B(\mathbb R).$
Conversely, if $f$ has a spectral factor $s \in B(\mathbb R)$ then 
$h := \chi_{b(f)/4}\, s \in H(\mathbb R)$ hence $P[h] = e^{u+iv}.$
Since $|h|^2 = f \geq m > 0,$ Lemma \ref{bohr3} implies that 
$v|_{\mathbb R} = \arg h$ has the form in Lemma \ref{bohr3}.
\\ \\
The following result shows that the assumption $f \in A(\mathbb R)$ in Lemma \ref{levin1} is not necessary to ensure that $f$ has a spectral factor $s \in B(\mathbb R).$
\begin{theo}\label{result2}
For $m > 0$ there exists $f \in B(\mathbb R)\, \backslash \, A(\mathbb R)$ with $f \geq m$ whose spectral factor $h \in B(\mathbb R)\, \backslash \, A(\mathbb R).$
\end{theo}
{\bf Proof} Rudin (\cite{rudin}, Theorem 5.12) gave a proof, based on the Banach-Steinhaus theorem
(\cite{rudin}, Theorem 5.8), that the subset of $U(\mathbb R)$ of period $2\pi$--periodic functions with non absolutely convergent Fourier series is nonempty. Zygmund (\cite{zygmund}, Chapter VI, 3.7)
gave the example
$$
	\phi(x) \sim \sum_{n = 2}^{\infty} \frac{\sin n x}{n \log n}.
$$
Fej\'{e}r (\cite{zygmund}, Chapter III, Theorem 3-4) proved its Ces\`{a}ro sums
$$
	p_n(x) = \sum_{k = 2}^n \frac{n+1-k}{n}\, \frac{\sin k x}{k \log k}, \ \ n \geq 2
$$
converge uniformly to $\phi$ therefore there exists an integer sequence $2 \leq n_1 < n_2 < \cdots$ with $||\phi-p_{n_j}|| \leq 2^{-j}/3.$  Define $q_1 = p_{n_1}$ and
$q_j = p_{n_{j+1}} - p_{n_j}, j \geq 2.$ Then $||q_j|| \leq 2^{-j}$ and $\min \Omega(q_j) = -n_{j+1}$ and $\max \Omega(q_j) = n_{j+1}.$ Construct a sequence $\rho_j  \in (0,n_{j+1}^{-1})$
with $\{\rho_j\}$ lineary independent over $\mathbb Q$ and define
$$
	g_j := \mathfrak{D}_{\rho_j} q_j, \ \ j \geq 1.
$$
Then $||g_j|| = ||q_j||, ||g_j||_A = ||q_j||_A,$ and $\Omega(g_j) = \rho_j \Omega(q_j)$
are pairwise disjoint subsets of $(-1,1).$ Therefore
$$
	g := \sum_{j = 1}^{\infty} g_j .
$$
satisfies $g \in B(\mathbb R), \Omega(g) \subset (-1,1),$ and
$$
||g||_A = \lim_{k \rightarrow \infty} \sum_{j = 1}^{k-1} ||q_j||_A \geq
\lim_{k \rightarrow \infty} ||\sum_{j = 1}^{k-1} q_j||_A =
\lim_{k \rightarrow \infty} ||p_{n_k}||_A = \infty
$$
so $g \in B(\mathbb R)\, \backslash \, A(\mathbb R),$ $g$ is real valued, and $\Omega(g) \subset (-1,1).$
Define $\Delta = -\inf \Omega(g),$ $h_1 := \chi_{\Delta}g,$ and $H_1 = P[h_1].$ 
Lemma \ref{asym} implies that $H_1$ is bounded and $h_1 \in B(\mathbb R) \cap H(\mathbb R).$
Define $c := \sqrt m - \inf \Re H,$ $h := c + h_1,$ $f := |h|^2,$ and $H := P[h].$ 
Lemma \ref{poisson} implies that $|H(z)| \geq |\Re H(z)| \geq \sqrt m$ so $H$ has no zeros in $\mathbb U.$
Then $s := \chi_{-\Delta} h$ is a spectral factor of $f.$ Since $s \notin A(\mathbb R),$ 
Lemma \ref{levin1} implies that $f \notin A(\mathbb R).$
\section{Research Questions}
We suggest the following questions for future research:
\begin{enumerate}
\item If the hypothesis $f \geq m > 0$ in Lemma \ref{levin1} and Theorem \ref{result1} is replaced by
the weaker hypothesis $f \geq 0,$ what conclusions about the spectral factor $s$ of $f$ can be deduced?
\item If $f \in B(\mathbb R)$ is nozero and $f \geq 0,$ Bohr showed that it lifts to a function $\widetilde f \in C(\mathbb R_B),$ where $\mathbb R_B$ is the Bohr compactification of $\mathbb R,$  and we proved in \cite{lawton} that $\log \widetilde f \in L^1(\mathbb R).$ Helson and Lowdenslager \cite{helsonlow} proved that  
$\widetilde f = |\widetilde h|^2$ where $\widetilde h$ is an outer function in the Hardy space $H^2(\mathbb R_B)$ (with respect to the linear order on the Pontryagin dual of $\mathbb R_B$ which equals the discrete real group). What is the relationship between $\widetilde h$ and the lift
$\widetilde s$ of the spectral factor $s$ of $f$ give by Lemma \ref{krein}?
\item Ahiezer's result in Lemma \ref{ahiezer} holds for operator valued and matrix valued functions \cite{ephremidze1,ephremidze2,gohberg1,gohberg2,rosenblumrovnyak}. What analogues do the results in this paper have in this context?
\end{enumerate}
{\bf Acknowledgments} The author thanks August Tsikh and Lasha Ephremidze for
discussions about entire functions and matrix spectral factorization.

\end{document}